\documentclass[12pt]{article}
\usepackage{cmap}
\usepackage[T2A]{fontenc}
\usepackage[cp866]{inputenc}
\usepackage[english,russian]{babel}
\usepackage[tbtags]{amsmath}
\usepackage{amsfonts,amssymb}

\newcommand\np{\noindent}

\newcommand\CC{\mathbb C}
\newcommand\DD{\mathbb D}

\newcommand\ZZ{\mathbb Z}
\newcommand\RR{\mathbb R}

\newcommand\ga{\gamma}
\newcommand\sub{\subset}
\newcommand\sm{\setminus}
\newcommand\pa{\partial}

\newcommand\al{\alpha}
\newcommand\ex{\ \exists\ }
\newcommand\la{\lambda}

\newcommand\eps{\varepsilon}

\newcommand\trl{\blacktriangleleft}
\newcommand\trr{\blacktriangleright}

\begin{document}

\title{\Large Harnack inequalities, Kobayashi distances\\
and holomorphic motions}

\author{E.M.Chirka}
%\thanks{The paper is supported by RFBR, gr. 11-01-00495
%and 11-01-00495-ofi-m.}}

\date{1 April 2012}

\maketitle

\np{\footnotesize  We prove some generalizations and analogies
of Harnack inequalities for pluriharmonic, holomorphic and
``almost holomorphic'' functions.
The results are applied to the proving of smoothness properties of
holomorphic motions over almost complex manifolds.}
%\end{abstract}
\vskip6mm

The Harnack inequalities in this paper are considered only from the point
of view of Complex Analysis, in application to real parts of holomorphic
functions.
In this context they can be written in invariant form independent of
biholomorphic transformations and thus can be generalized onto arbitrary
almost complex manifolds what we show in sec.1.

First of all we are interesting in holomorphic functions which do not
take the values 0 and 1 because just such functions constitute in the main
the normalized holomorphic motions (see sec.$\,$3).
The main estimates of such functions give the theorems of Landau and
Schottky which were specified many times and are presented in the final form
in [H].
The estimates of J.Hempel can also be written in invariant form and thus
they are evidently extendable onto general almost complex manifolds
(sec.$\,$2).

Most likely, these results are already known but the proofs are so simple
that it seemed for me more complicated to look for precise references, and
I do not affirm that these results are new.

The general Landau theorem is applied in sec.$\,$3
in the proving of an analogy
of Harnack inequalities for holomorphic functions with values in
$\CC\sm\{0,1\}$ (sec.$\,$2) and (following [GJW]) in the
proof of H\"older conditions for arbitrary holomorphic motions (sec.$\,$3).

In sec.$\,$4 we establish analogies of Harnack inequalities
for functions in disk
or on the plane which satisfy special estimates of derivatives in
$\bar z$ what has direct relation with almost holomorphic motions, in
particular, with holomorphic motions in nonstandard complex structures.
\vskip5mm

{\bf 1.$\,$Pluriharmonic functions.$\,$}
Let us start with base Harnack inequalities for positive harmonic functions
$u(z)$ in the unit disk $\DD:|z|<1$ on the complex plane $\CC\,$:
\begin{equation}
\frac{1-|z|}{1+|z|}\,\le\,\frac{u(z)}{u(0)}\,\le\,\frac{1+|z|}{1-|z|}\,
\label{1}
\end{equation}
(see e.g. [HK]).
Noticing that the Poincar\'e distance $\,\rho_{\,\DD}$ between the points
$0\,$ and $\,z$ in $\,\DD$ is equal to $\,log\,\frac{1+|z|}{1-|z|}$ (assuming
that Poincar\'e metric has curvature
${}\equiv -1$) we can rewrite these inequalities in invariant form independent
of M\"obius transformations of $\,\DD\,$:
$$
e^{-\rho_{\,\DD}(z,z_0)}\,\le\,\frac{u(z)}{u(z_0)}\,
\le\,e^{\,\rho_{\,\DD}(z,z_0)}
$$
for any $z,z_0\in\DD$.
\vskip2mm

But the metric of curvature ${}\equiv -1$ exists on arbitrary hyperbolic
Riemann surface, and it is natural to ask if the correspondent inequalities
are valid on such surfaces.
And if one notices that harmonic functions in the disk are partial case of
pluriharmonic functions on a complex manifold then it is natural to interest
in analogy of the last inequalities for such functions also.
\vskip1mm

A generalization of Poincar\'e distance on arbitrary complex manifold is the
Kobayashi distance (see [K]) wich coincides in the disk with $\,\rho_{\,\DD}$.
Consider it more detaily in very general setting.
\vskip2mm

Let $B$ be arbitrary path connected complex Banach manifold
(in particular, finite dimensional one).
For any $z,z'$ placed in one coordinate ball of the manifold
$B$ there exists evidently a holomorphic disk $h:\DD\to B$ such that
$h(0)=z',\,h(\zeta)=z$ for some $\zeta\in\DD$.
As $B$ is path connected and any path $[0,1]\to B$ is covered by a finite number
of coordinate balls, then for any $z,z'\in B$ there exists a chain of holomorphic
disks $h_j:\DD\to B$ and points $\zeta_j\in\DD$ such that
$h_1(0)=z',\,h_{j+1}(0)=h_j(\zeta_j),\,j=1,...,N-1
$, and $h_N(\zeta_N)=z$.
The amount
$$
\kappa_B(z,z')=\inf\ \sum_1^N\,log\,\frac{1+|\zeta_j|}{1-|\zeta_j|}
$$
where infimum is taken by all described chains of holomorphic disks
($N$ can be arbitrary) is called {\it\/Kobayashi distance\/} between the
points $z,z'$ on the manifold $B$.
As $\,log\,\frac{1+|\zeta|}{1-|\zeta|}=2\,arcth\,|\zeta|$
is the Poincar\'e distance in $\DD$ between $0,\zeta$ and this distance is
invariant with respect to M\"obius automorphisms of $\DD$ then
$\kappa_B(z,z')=\kappa_B(z',z)$.
\vskip2mm

In general case $\kappa_B$ is only a pseudometric, the distances between some
different points can be equal to zero but the triangle inequality follows easily
from the definition.

The most important property of Kobayashi distance is its evident
non\-increasing
by holomorphic mappings: if $f:B\to X$ is such a map then
\begin{equation}
\kappa_{X}(f(z),f(z'))\,\le\,\kappa_B(z,z')
\label{2}
\end{equation}

\np (it is one of abstract variants of {\it\/Schwarz lemma\/}).
\vskip3mm

Kobayashi distance can be defined as above on arbitrary almost complex Banach
manifold $B$.
As the existence of finite chains of holomorphic disks connecting given
$z,z'\in B$ on such a manifold is not evident at all, then we set for
definitness that $\kappa_B(z,z')=\infty$ if there is no such a chain.
For marked (base) point $z_0$, the ball $\{z\in B:\kappa_B(z,z_0)<R\}$
will be denoted by $B_R\,$.
\vskip3mm

Let us remind that a continuous function $u:B\to\RR$ on almost complex Banach
manifold is called {\it\/pluriharmonic\/} if for any holomorphic map
$f:\DD\to B$ the function $u\circ f$ is harmonic in the unit disk.
\vskip3mm

{\bf Proposition 1.} {\it Let $B$ be an almost complex Banach manifold and
$\,u$ be a positive pluriharmonic function on $B$.
Then, for any} $z,z_0\in B$,
\begin{equation}
e^{-\kappa_B(z,z_0)}\,\le\,\frac{u(z)}{u(z_0)}\,
\le\,e^{\,\kappa_B(z,z_0)}\,.
\label{3}
\end{equation}
\vskip1mm

$\trl$ If
$\kappa_B(z,z_0)=\infty$ then there is nothing to prove, thus we assume further
that $\kappa_B(z,z_0)<\infty$.

Let us fix $\eps>0$ and choose a chain of holomorphic disks $h_j:\DD\to B$ and
$\zeta_j\in\DD$ such that
$h_1(0)=z_0,\,h_{j+1}(0)=h_j(\zeta_j),\,j=1,...,N-1,\,h_N(\zeta_N)=z$ and
$\Sigma_1^N\,log\,\frac{1+|\zeta_j|}{1-|\zeta_j|}\le\kappa_B(z,z_0)+\eps$.

By the classical Harnack inequality
$$
\frac{1-|\zeta_j|}{1+|\zeta_j|}\le
\frac{u\circ h_j(\zeta_j)}{u\circ h_j(0)}\le\frac{1+|\zeta_j|}{1-|\zeta_j|}\,.
$$
As
$$
\frac{u(z)}{u(z_0)}=\frac{u\circ h_N(\zeta_N)}{u\circ h_N(0)}
\cdots\frac{u\circ h_1(\zeta_1)}{u\circ h_1(0)}\,,
$$
it follows that
$$
e^{-\kappa_B(z,z_0)-\eps}\le\frac{u(z)}{u(z_0)}
\le e^{\,\kappa_B(z,z_0)+\eps}
$$
and the inequalities are proved due to arbitrary $\eps$.
$\trr$
\vskip2mm

The inequalities (3) are applicable to the functions of type
$\,log\,M/|f|$ where $f$ is a holomorphic function with $0<|f|<M$.
Consider as an example a ``pointwise'' analogy of two constant theorem for such
functions in which the point $z_0$ plays the role of a set of positive harmonic
measure.
\vskip2mm

{\bf Corollary.} {\it If $f$ is a holomorphic function on $B$ and $0<|f|<M$ then
$$
|f(z)|\,<\,|f(z_0)|^{\,\al(z)} M^{1-\al(z)},\ \,z\in B
$$
where $\al(z)=e^{-\kappa_B(z,z_0)}$.}
\vskip2mm

$\trl$
The function $log\,M/|f|$ is pluriharmonic and positive.
According (3), $M/|f(z)|\ge(M/|f(z_0)|)^{\,\al(z)}$ what is clamed.
$\trr$
\vskip5mm

{\bf 2.$\,$Functions with values in $\CC\sm\{0,1\}$.}
The hyperbolic domain $\CC\sm\{0,1\}$
with complete hyperbolic Poincar\'e metric $\rho_{\,0,1}$
plays big role in different problems of Complex Analysis.
By the universal covering $\DD\to\CC\sm\{0,1\}$ this $\,\rho_{\,0,1}$
is lifted to Poincar\'e metric $\rho_{\,\DD}$.
As the lifting to $\DD$ of Kobayashi metric on $\CC\sm\{0,1\}$
is the same then
$\rho_{\,0,1}$ {\sl\/coincides\/} with Kobayashi metric in this domain.
\vskip2mm

Infinitesimal Poincar\'e metric in $\CC\sm\{0,1\}$ has the form
$\rho_{\,0,1}(z)|dz|$ (i.e.
$\rho_{\,0,1}(z_1,z_2)=\inf_\ga\int_\ga\rho_{\,0,1}(z)|dz|$
for any $z_1,z_2$
where infimum is taken by all smooth pathes with ends $z_1,z_2$).
We will need in the following estimates of this metric.
\vskip3mm

{\bf Lemma 1.} {\it $\ \,\rho_{\,0,1}(z)\ge\rho_{\,0,1}(-|z|)\ge
\left(|z|(C_{0,1}+\,log\,\frac{1}{|z|})\right)^{-1}$
in $\,\DD\sm 0\,$ with the constant $\,C_{0,1}=1/\rho_{\,0,1}(-1)>1$.}
\vskip2mm

$\trl$
The first inequality is proved in the paper of Lehto and Virtanen [LV],$\,$p.6
(see also [Ag]).

The second inequality is proved in the paper of Hempel [H],$\,$p.443, but we
present here a simpler proof.

It is wellknown (see e.g. [A1],$\,$1-8) that
$$
log\,\rho_{\,0,1}(z)+log\,\frac{1}{|z|}+log\,(log\frac{1}{|z|})\to 0
$$
as $z\to 0$.
(More detailed asymptotics see in [H].)

As $min_{\,\pa\bar\DD}\,\rho_{\,0,1}=\rho_{\,0,1}(-1)$ then
$\rho_{\,0,1}(z)\ge \rho_{\,0,1}(-1)=1/(|z|\,log\,\frac{C}{|z|})$ when $|z|=1$
with the constant $C=exp\,1/\rho_{\,0,1}(-1)$.
Setting $\rho_0(z)=1/(|z|\,log\,\frac{C}{|z|})$ we obtain that
$log\,\rho_0(z)-log\,(1/(|z|\,log\,\frac{1}{|z|}))\to 0$ as $z\to 0$, hence,
due to the asymptotic obtained above,
$log\,\rho_{\,0,1}(z)-\,log\,\rho_0(z)\to 0$ as $z\to 0$.

As the metrics $\rho_{\,0,1}(z)$ and $\rho_0(z)$ have the same Gauss curvature
($\Delta\,log\,\rho\equiv\rho^2$) then
$\Delta\,(\,log\,\rho_{0,1}-log\,\rho_0)=\rho_{0,1}^2-\rho_0^2$.
The function $\rho_{\,0,1}/\rho_0\ge 1$ on $\pa\DD$ by the definition of $C$, and
the same at $0$ in sense of limit as it is proved above.
If it takes minimum at some point $z_0\in\DD\sm 0$ then there will be also
the miminum of $log\,\rho_{\,0,1}-log\,\rho_0$.
At the minimum point of a smooth function its laplasian is non-negative, hence
$\rho_{\,0,1}^2(z_0)\ge\rho_0^2(z_0)$ and thus $\rho_{\,0,1}\ge\rho_0$ everywhere
in $\DD\sm 0$.

Comparing $\rho_{\,0,1}$ with Poincar\'e metric for half-plane
$\{Re\,z<0\}\sub\CC\sm\{0,1\}$ we obtain that $C_{\,0,1}>1$.
$\trr$
\vskip3mm

Another proof of this lemma (with a bigger constant $C$) is contained in [GJW]
where we have taken the idea of using Gauss curvatures going back to
Ahlfors (see [A1],$\,$1-5).

By Agard formula [Ag]
$$
\frac{1}{\rho_{\,0,1}(z)}=
\frac{1}{2\pi}\int_\CC \left|\frac{z(z-1)}{\zeta(\zeta-1)(\zeta-z)}\right|
dS_\zeta
$$
and thus the precise value of the constant $C_{0,1}$ is equal to
$\frac{1}{\pi}\int_\CC|\zeta(\zeta^2-1)|^{-1}dS_\zeta$.
A simpler expression
\vskip1mm

\centerline{$C_{\,0,1}=\Gamma(\frac{1}{4})^4/4\pi^2=4,3768796...$}
\vskip2mm

\np is given in [H], a proof can be found in [N],$\,$Ch.6,$\,$sec.6.
\vskip3mm

Kobayashi metric is invariant with respect to holomorphic automorphisms, in
particular, the metric $\,\rho_{\,0,1}(z)|dz|$ is invariant with
respect to the transfor\-mation $z\mapsto 1/z$.
From this and infinitesimal Schwarz lemma (see (5) below) for imbeddings of
$\DD\sm 0$ and $\CC\sm\bar\DD\,$ into $\,\CC\sm\{0,1\}$ we obtain
\vskip3mm

{\bf Corollary 1.}
{\it For all $z\in\CC\sm\{0,1\},$}
\begin{equation}
|z|\,|\,log\,|z||\le\frac{1}{\rho_{\,0,1}(z)}
\le |z|(\,C_{\,0,1}+|\,log\,|z||)\,.
\label{4}
\end{equation}

The same estimates are valid also for Poincar\'e metrics $\,\rho_{\,0,a}$ in
$\CC\sm\{0,a\}$ with arbitrary
$\,a,\,|a|=1$, because $\rho_{\,0,a}(z)=\rho_{\,0,1}(z/a)$ and
$C_{\,0,1}=1/\rho_{\,0,a}(-a)$.
\vskip4mm

Infinitesimal form of Kobayashi distance on arbitrary finite dimensional
complex manifold $B$ is the {\it\/Royden metric\/} (or Kobayashi -- Royden)
which mesures the lengths of tangent vectors $V\in T_z B\,$ to $B$
by the formula
\vskip2mm

\centerline{{\it $|V|_\kappa\equiv |V|_{\kappa(B)}=inf\,\{2/R:\ex$
holomorphic disk
$\,h:\DD_\zeta\to B,\,\zeta=\xi+i\eta$,}}
\centerline{{\it such that
$\,h(0)=z\,$ and $\,h_*(\frac{\pa}{\pa\xi}|_0)=R\,V\}$,}}
\vskip2mm

\np where $h_*:T\DD\to TB$ is the tangent map to $h$.
\np (If Poincar\'e metric in the disk is normalized in another way, as
$\frac{|d\zeta|}{1-|\zeta|^2}$ with curvature ${}\equiv -4$,
then instead of $2/R$ in the definition one has to put $1/R$).
In particular, in the disk $\DD$ the Kobayashi -- Royden metric
coincides with the Poincar\'e metric, $|V|_{\kappa(\DD)}=2\,|V|/(1-|z|^2)$
for all $V\in T_z\DD$.

It is, as a rule, not Riemannian but Finsler (semi)metric (there can be
nonzero vectors of zero length) and
$\kappa_B(z,z')=\inf\,\int_0^1|\ga'(t)|_\kappa\,dt\,$
where infimum is taken by all smooth pathes $\ga:[0,1]\to B$ with ends
$z,z'$ (see [R1,R2]).
\vskip3mm

The definition of Kobayashi -- Royden metric given above suits surely for
arbitrary almost complex Banach manifolds but I do not know if the
integral representation given above is valid for almost complex and infinite
dimensional manifolds.
Nevertheless, also in general case we have
{\it\/infinitesimal Schwarz lemma\/}:
\begin{equation}
|f_*(V)|_{\kappa(Y)}\,\le\,|V|_{\kappa(X)}\,,\quad V\in TX\,,
\label{5}
\end{equation}
for any holomorphic map $f:X\to Y$ of almost complex Banach manifolds.
The proof follows evidently from the definition of Kobayashi -- Royden metric.

For example, if $B$ is the ball $\|z\|<R$ in complex Banach space
and $h:\DD\to B,\,h(0)=0,$ is a holomorphic disk then
$\|h(\zeta)\|\le R\,|\zeta|$ by Schwarz lemma, hence $\|h_*(0)\|\le R$
($h_*(0):TB\to X$ because $TB=B\times X$).
If $h_*(\frac{\pa}{\pa\xi}|_0)=R\,'\,V$ then it follows that
$R\,'\,\|V\|\le R,\,2/R\,'\ge 2\|V\|/R$, hence
$|V|_{\kappa(B)}\ge 2\,\|V\|/R$.
The extremal disk $h:\DD\to B$ realising the infimum
is the linear map $\zeta\mapsto R\,\zeta\,V/\|V\|$;
therefore $|V|_\kappa=2\,\|V\|/R$ for all $V\in T_0B$.
\bigskip

We apply these notions to holomorphic functions with values in
$\CC\sm\{0,1\}$.
\vskip3mm

{\bf Proposition 2.} {\it Let $f$ be a holomorphic function on almost complex
Banach manifold $B$ with values in $\CC\sm\{0,1\}$ and $\,V\in T_z B\,$.
Then
$$
|(Vf)(z)|\le |V|_\kappa\cdot |f(z)|\,(\,C_{0,1}+|\,log\,|f(z)||)
$$
where $\,C_{\,0,1}=1/\rho_{0,1}(-1)>1$.}
\vskip2mm

$\trl$
Let us fix arbitraty $\eps>0$.
Let $h:\DD\to B$ be a holomorphic disk such that
$\,h(0)=z,\,h_*(\pa/\pa\xi|_0)=2\,V/(|V|_\kappa+\delta)$ with
$0\le\delta\le\eps$.
Then $g:=f\circ h:\DD\to\CC\sm\{0,1\}$.
By infinitesimal Schwarz lemma, for any such map $g$ the following
inequality is valid:
$$
\rho_{0,1}(g(\zeta))\,|dg(\zeta)|\le 2\,|d\zeta|/(1-|\zeta|^2)\,,
$$
in particular, $\,|g'(0)|\le 2/\rho_{0,1}(g(0))$.
As $g(0)=f(z)$ then by Lemma 1
$$
|g'(0)|\le 2\,|g(0)|\,(\,C_{\,0,1}+|\,log\,|g(0)||)=
2\,|f(z)|\,(\,C_{\,0,1}+|\,log\,|f(z)||).
$$
By the choice of $h$,
$$
g'(0)=\left.\frac{\pa}{\pa\zeta}\right|_0 (f\circ h)(\zeta)=
\left.\frac{\pa}{\pa\xi}\right|_0 (f\circ h)(\zeta)=
\frac{2(Vf)(z)}{|V|_\kappa+\delta}\,,
$$
hence $|(Vf)(z)|\le(|V|_\kappa+\eps)\cdot
|f(z)|\,(\,C_{\,0,1}+|\,log\,|f(z)||)$
and the statement is proved due to arbitrary $\eps$.
$\trr$
\vskip3mm

If $B=\DD\,$ and $\,V=\frac{\pa}{\pa\xi}|_0\,$
then $\,|V|_\kappa=\frac{2}{1-|z|^2}\,$, hence
$$
|f'(z)|\le\frac{2}{1-|z|^2}\cdot |f(z)|\,
(\,C_{0,1}+|\,log|f(z)||)\,.
$$
By $z=0$ it is classical Landau theorem:
\vskip2mm

{\bf Corollary 2.} {\it Let $f(z)=\Sigma_{\,0}^{\infty}\,a_n\,z^n$
be a function
holomorphic in $\DD$ which does not take the values $\,0$ and $1$.
Then}
$$
|a_1|\le 2\,|a_0|\,(C_{\,0,1}+|\,log\,|a_0||)\,.
$$

The exact constant $\,C_{\,0,1}=1/\rho_{0,1}(-1)\,$ in this theorem is
established by J.Hempel [H]; the equality is attained by the universal
covering $\,S:\DD\to\CC\sm\{0,1\},\,S(0)=-1$.
\vskip2mm

Due to the estimate of Kobayashi -- Royden metric in a ball we obtain
similar estimate.
\vskip2mm

{\bf Corollary 3.} {\it Let $f$ be a holomorphic function in the unit ball
of a Banach space $(X,\|\cdot\|)$ which does not take the values $0,1$.
Then}
$$
\|f'(0)\|\le 2\,|f(0)|\,(C_{\,0,1}+|\,log\,|f(0)||)\,.
$$
\vskip2mm

Thus, Proposition 2 is a generalization of Landau theorem onto arbitrary
almost complex manifolds and, as we see, the proof is practically the same
as for $\,B=\DD$ and is a simple corollary of infinitesimal Schwarz lemma
(5) and the estimate (4) of Poincar\'e metric $\,\rho_{\,0,1}$.
\vskip3mm

The following are analogies of Harnack inequalities (3).
\vskip2mm

{\bf Proposition 3.} {\it
Let $f$ be a holomorphic function on almost complex Banach manifold $B$
with values in $\CC\sm\{0,1\}$.
Then
\begin{equation}
e^{-\kappa_B(z,z_0)}\le
\frac{C_{\,0,1}+|\,log\,|f(z)||}{C_{\,0,1}+|\,log\,|f(z_0)||}\le
e^{\,\kappa_B(z,z_0)}
\label{6}
\end{equation}
for any $z,z_0\in B$.
Furthermore, if there exists continuous $\,log\,f$ on $B$ then}
\begin{equation}
e^{-\kappa_B(z,z_0)}\le
\frac{C_{\,0,1}+|\,log\,f(z)|}{C_{\,0,1}+|\,log\,f(z_0)|}\le
e^{\,\kappa_B(z,z_0)}
\label{7}
\end{equation}
\vskip2mm

$\trl$
Let us consider first the crucial case $B=\DD\sub\CC_\zeta\,$ with
$\,\zeta=re^{it}$; here continuous logarithm exists.

Fix an arbitrary point $a=|a|\,e^{i\al}\in\DD$ and denote
$u(r):=|\,log\,|f(re^{i\al})||$, $C:=C_{\,0,1}$.
If $u$ is not constant (what we assume further) then there can be on
the interval $(0,a)$ only finite set of points where
$u=0$, hence the function $u$ is piecewise smooth.
As
$$
\left|\frac{\pa}{\pa r}\,log\,f\right|
\ge\left|\frac{\pa}{\pa r}\,log\,|f|\right|\ge\frac{\pa}{\pa r}|\,log\,|f||
$$
almost everywhere on $(0,a)$ then by Prop.$\,$2
$$
\frac{u'}{C+u}\le\frac{2}{1-r^2}=\frac{1}{1+r}+\frac{1}{1-r}\,.
$$
Integrating this inequality by $(0,|a|)$ we obtain that
$$
log\,\frac{C+u(|a|)}{C+u(0)}\le\,log\,\frac{1+|a|}{1-|a|}\,.
$$
The substitution $\zeta=\frac{a-\eta}{1-\bar a\eta}$ changes $0,a$ in places,
hence
$$
\frac{1-|a|}{1+|a|}\le\frac{C+|log\,|f(a)||}{C+|log\,|f(0)||}\le
\frac{1+|a|}{1-|a|}\,.
$$

For obtaining inequalities (7) let us notice that there exists continuous
$\,log\,(\,log\,f)$ on $[0,a]$.
As $|\frac{\pa}{\pa r}\,log\,(\,log\,f)|\ge\frac{\pa}{\pa r}\,log|\,log\,f|$
then
$|\frac{\pa}{\pa r}log\,f|\ge\frac{\pa}{\pa r}|\,log\,f|$.
Thus setting $v(r)=|\,log\,f(re^{i\al})|$ we obtain, again by Prop.$\,$2,
that $v'/(C+v)\le 2/(1-r^2)$.
Integrating this as above we obtain the inequalities (7) in $\DD$.
\vskip1mm

In the case of arbitrary base $B$ we repeat the argument from the proof
of Prop.$\,$1.
If $\kappa_B(z,z_0)=\infty$ then there is nothing to prove, therefore we
assume further that $\kappa_B(z,z_0)<\infty$.

Fix $\eps>0$, choose a chain of holomorphic disks
$h_j:\DD\to B$ and points $\zeta_j\in\DD$, $z_0,z_1,...,z_N=z$ such that
$h_j(0)=z_{j-1},\,h_j(\zeta_j)=z_j,\,j=1,...,N$, and
$\,\kappa_B(z,z_0)\ge\Sigma_1^N\,log\,\frac{1+|\zeta_j|}{1-|\zeta_j|}-\eps\,.$
As it is proved above,
$$
\frac{1-|\zeta_j|}{1+|\zeta_j|}\le
\frac{C+|log\,|f(z_j)||}{C+|log\,|f(z_{j-1})||}\le
\frac{1+|\zeta_j|}{1-|\zeta_j|}
$$
and corresponding inequalities for $|\,log\,f|$ if there exists on $B$
a continuous logarithm of $f$.

It follows that
$$
e^{-\kappa_B(z,z_0)-\eps}\le
\frac{C+|\,log\,|f(z)||}{C+|\,log\,|f(z_0)||}=
\prod_1^N\frac{C+|\,log\,|f(z_j)||}{C+|\,log\,|f(z_{j-1})||}\le
e^{\,\kappa_B(z,z_0)+\eps}
$$
with arbitrary $\eps>0$, and corresponding inequalities for $\,|log\,f|$ if
a continuous logarithm of $f$ exists.
$\trr$
\vskip3mm

Using for metrics $\rho_{\,0,a}$ with $|a|=1$ the same estimates (4) as for
$\rho_{\,0,1}$ (see above) we obtain similarly the inequalities
$$
e^{-\kappa_B(z,z_0)}\le
\frac{C_{\,0,1}+|\,log\,e^{-i\,arg\,f(z_0)}f(z)|}
{C_{\,0,1}+|\,log\,|f(z_0)||}\le
e^{\,\kappa_B(z,z_0)}\,
$$
which give more precise estimate of the argument of
the function $f$.
\vskip2mm

The constant $C_{\,0,1}$ in the estimates (6) and the last inequalities
is exact.
Indeed, let $S:\DD\to\CC\sm\{0,a\},\,|a|=1,$ be the universal covering
normalized by conditions $S(0)=-a,\,S'(0)>0$.
Then $\rho_{\,0,a}(S(z))|S'(z)|=2/(1-|z|^2)$, hence
$\rho_{\,0,a}(-a)S'(0)=2$.
As $\,C_{\,0,1}=1/\rho_{\,0,a}(-a)=S'(0)/2$ then the right inequality can be
rewritten in the form
\vskip1mm

\centerline{$1+|\,log\,(1-S'(0)\,z\,(1+o(1)))|/C_{\,0,1}=
1+2\,|z|\,(1+o(1))\le
\frac{1+|z|}{1-|z|}\,,$}
\vskip1mm

\np and it is arbitrary close to equality for small $|z|$.
\vskip2mm

{\bf Corollary 4.} {\it Let $h$ be a holomorphic function on $B$ which does
not take values from $\,2\pi i\,\ZZ\,$.
Then}
$$
|h(z)|\le(\,C_{\,0,1}+|\,Re\,h(z_0)|)\cdot e^{\,\kappa_B(z,z_0)}+
|\,Im\,h(z_0)|-C_{\,0,1}\,.
$$

$\trl$
Let $\,f=e^h$.
Then $h=log\,f,\ log\,|f|=Re\,h$ and $\,arg\,f=Im\,h$.
Put this in the last right inequality.
$\trr$
\vskip2mm

Similar estimates are valid surely for functions which do not take values
in an arithmetic progression $\,a+b\,\ZZ,\,b\not=0$.
Therefore the growth of such a function $f(z)$, say in the unit ball, does
not exceed $\,C/(1-\|z\|)$ with corresponding constant $C=C(f,a,b)$.
\vskip3mm

Proposition 3 permits to estimate the growth of a function
$B\to\CC\sm\{0,1\}$ by its value in a fixed point and the distance to this
point in Kobayashi metric.
In the case $B=\DD$ the classical Schottky' theorem states that, for such
functions, $|f(z)|$ is estimated by a quantity depending only on $|f(0)|$
and $|z|$.
In general, let us denote by $M(R,R\,')$ the supremum of numbers $|f(z)|$
by all holomorphic functions $f:B\to\CC\sm\{0,1\}$ such that
$|f(z_0)|\le R\,'$ and by all $z\in B_R$.
\vskip3mm

{\bf Corollary 5.} {\it If $|f(z_0)|\le R\,'$ and
$\,z\in B_R\,$ then}
$$
|f(z)|\,\le\,M(R,R\,')\le e^{-C_{\,0,1}}(e^{\,C_{\,0,1}}(max(1,R\,'))^{\,e^R}\,.
$$
\vskip2mm

$\trl$
If $|f(z_0)|\ge 1$ then by Prop.$\,$2
$$
|\,log\,|f(z)||\le e^R(C_{\,0,1}+log\,R\,')-C_{\,0,1}\,,
$$
hence
$$
|f(z)|\le exp\,((C_{\,0,1}+log\,R\,')\,e^R-C_{0,1})=
e^{-C_{\,0,1}}(e^{\,C_{\,0,1}}\,R\,')\,e^R\,.
$$
It remains to notice that $M(R,R\,'')\le M(R,R\,')$ if $R\,''<R\,'$ and
thus $M(R,1)\le e^{\,C_{\,0,1}(e^R-1)}$.
$\trr$
\vskip2mm

This estimate is far from being exact even for the case $f(z_0)=-1$.
An exact but implicit estimate for $B=\DD$ is pointed by Hempel [H]$\,$:
$$
\left|\int_{|f(0)|}^{|f(z)|}\,\rho_{\,0,1}(-r)\,dr\right|\le
\frac{1+|z|}{1-|z|}\,.
$$
It follows evidently from the inequality
$\,\rho_{\,0,1}(z)\ge\rho_{\,0,1}(-|z|)$ (see [LV],[Ag]) and Schwarz lemma (2)
and thus it is generalized to holomorphic functions $f:B\to\CC\sm\{0,1\}$
on arbitrary almost complex Banach manifold$\,$:
$$
\left|\int_{|f(z_0)|}^{|f(z)|}\,\rho_{\,0,1}(-r)\,dr\right|\le
\kappa_B(z,z_0)\,.
$$

The equality here with $B=\DD$ is attained by any universal covering
$F:\DD\to\CC\sm\{0,1\}$ such that $F(0)<0$ and $F(z)<0$.
\vskip2mm

To estimate $M(R,R\,')$ by this formula seems rather difficult.
Several simple explicit estimates for $B=\DD$, such as that one in
Corollary 5, are contained in the paper of Jenkins [J] but I do not see
how to extend them to the case of arbitrary $B$.
\vskip4mm

At the end of the section,
several words on functions with values in $\DD\sm 0$.
Instead of estimates (4) we have here precise formula for Poincar\'e
metric,
$\rho_{\,\DD\sm 0}(z)=1/|z|\,log\,1/|z|$.
Repeating the proof of Prop.$\,$2 we obtain inequalities
$$
e^{-\kappa_B(z,z_0)}\le
\frac{log\,1/|f(z_0)|}{log\,1/|f(z_0)|}\le
e^{\,\kappa_B(z,z_0)}
$$
for holomorphic functions $f:B\to\DD\sm 0$ on almost complex Banach
manifold $B$ as well as other similar inequalities when there exists
a continuous logarithm of $f$.
\vskip5mm

{\bf 3.$\,$On the smoothness of holomorphic motions.}
{\it Holomorphic motion\/} of a set $E\sub\hat\CC$ over almost complex Banach
manifold $B$ with a base point $z_0$ is an arbitrary map
$\phi:B\times E\to\hat\CC$ with the following properties:

(1) $\,\phi(z_0,w)\equiv w\,$,

(2) $\,\phi(z,\cdot):E\to\hat\CC$ is injective for every fixed $z\in B\,$,

(3) $\,\phi(\cdot,w)\,$ is holomorphic in $B$ for every fixed $w\in E\,$.

\np A motion $\phi$ is called {\it normalized\/} if
$w=0,1,\infty\in E$ and $\phi(\la,w)\equiv w$ for these $w$.

\vskip2mm

The most important and basic for the proofs of properties of
holomorphic motions is the case $B=\DD,\,z_0=0$.
Just in this form holomorphic motions have appeared for the first time in
the paper of R.M\`a\~n\'e, P.Sad and D.Sullivan
[MSS] where it was proved so called $\la$-lemma (the term is already
excepted) which states that holomorphic motions are continuous by joint
variables and quasiconformal for every fixed $z\in B$.

Any quasiconformal mapping, say, of the plane $\CC$ onto itself can be
included into some holomorphic motion of the set $E=\CC$ as a map
$\phi(z,\cdot):\CC\to\CC$ with suitable $z\in\DD\,$.
As quasiconformal mappings do not satisfy in general Lipschitz condition
(even locally) one can not wait from holomorphic motions a smoothness in
$w$ better than H\"older one.
At the same time it was noticed in the paper [GJW]
that the estimates as in Landau theorem imply the H\"older conditions
for normalized holomorphic motions.
We expose this below in more general context, with some specifications.
\vskip3mm

If $f_1,f_2$ are different functions from correspondent normalized
holomorphic motion then
$f_1-f_2$ in general can take value 1 and by this reason the estimate
of ``derivatives'' of the function $f_1-f_2$ comparing with
$f_j$ itself is a little more complicated.
\vskip3mm

{\bf Lemma 2.} {\it Let $f_1,f_2$ be holomorphic functions on almost
complex Banach manifold with values in $\CC\sm\{0,1\}$ which are different
at every point of $B$.
Then for any vector field $\,V$ on $B$ the following inequality is valid}
$$
\frac{|V(f_1-f_2)|}{|f_1-f_2|}\le |V|_\kappa\cdot
(\,2\,C_{\,0,1}+2\,min\,(|\,log\,|f_1||,|\,log|1-f_1|)+|\,log\,|f_1-f_2||)\,.
$$
\vskip2mm

$\trl$
The function $1-f_2/f_1$ maps $B$ to $\CC\sm\{0,1\}$, hence Prop.$\,$2 is
applicable to it.
As
$$
V\left(\frac{f_1-f_2}{f_1}\right)\cdot\frac{f_1}{f_1-f_2}=
\frac{f_2\,Vf_1-f_1\,Vf_2}{f_1(f_1-f_2)}=
\frac{Vf_1-Vf_2}{f_1-f_2}-\frac{Vf_1}{f_1}
$$
and $|Vf_1|/|f_1|\le|V|_\kappa(C_{\,0,1}+|log|f_1||)$ by Prop.$\,$2 then
$$
\left|\frac{Vf_1-Vf_2}{f_1-f_2}\right|\le|V|_\kappa\left(\,2\,C_{\,0,1}+
\left|\,log\,\left|\frac{f_1-f_2}{f_1}\right|\right|+|\,log\,|f_1||\right)
$$
$$
\le\,|V_\kappa|(\,2\,C_{\,0,1}+2\,|\,log\,|f_1||+|\,log\,|f_1-f_2||)\,
$$
and it remains to notice that this inequality is valid also for functions
$1-f_1,\,1-f_2$ with the difference $f_2-f_1$.
$\trr$
\vskip3mm

Lemma 2 shows that, for any normalized holomorphic motion
$\phi:B\times E\to\CC$
and any Kobayashi bounded vector field $V$ on $B$, the function $V\phi$
has almost the same modulus of continuity in $w$ as the motion $\phi$ itself.
This permits to establish the following estimate which can not be
improved in general.
\vskip2mm

{\bf Proposition 4.} {\it Let $\phi:B\times E\to\CC$ be a normalized
holomorphic motion of a set $E\sub\CC$ over
almost complex Banach manifold $(B,z_0)$.
Then for any $w_1,w_2\in E$, any $R>0$ and $z\in B_R$ the following
inequality is valid $($H\"older inequality with respect to spherical metric
on $\hat\CC)\,:$
\begin{equation}
\frac{|\phi(z,w_1)-\phi(z,w_2)|}
{\sqrt{1+|\phi(z,w_1)|^2}\sqrt{1+|\phi(z,w_2)|^2}}\,
\le\,
C_R\left(\frac{|w_1-w_2|}{\sqrt{1+|w_1|^2}
\sqrt{1+|w_2|^2}}\right)^{e^{-R}}
\label{8}
\end{equation}
with constant $C_R$ depending only on $R\ ($and {\bf not} depending on
$\phi\,)$.
Moreover, for any $R\,'>0$ there exists a constant $C(R,R\,')>1$ such that
\begin{equation}
(|w_1-w_2|/C(R,R\,'))^{e^R}\le |\phi(z,w_1)-\phi(z,w_2)|
\le C(R,R\,')\,|w_1-w_2|^{e^{-R}}
\label{9}
\end{equation}
by the condition that $|w_1|,|w_2|<R\,'$.}
\vskip2mm

Let us note at once that for the unit ball $B$ in a complex Banach space
and $z_0=0$ the exponent
$e^{-R}$ over the ball $\|z\|<r$ is equal to $\frac{1-r}{1+r}$ and this
estimate can not be improved in general.
\vskip2mm

$\trl$
Fix $R>0$, denote $\,e^{-R}=\al,\,\phi(\cdot,w_j)=f_j,\,f=f_1-f_2\,$
and prove at first the inequality (8).

If $|w_1|\le 1/2$ and $|w_2|>2$ or vice versa
then $|w_1-w_2|/\sqrt{1+|w_1|^2}\sqrt{1+|w_2|^2}>1/2$ and inequality
(8) is valid with any constant $C_R>2$ as its left side ${}\le 1$.
Therefore we can assume that both $|w_j|\le 2$ or both $|w_j|\ge 1/2$.
As the inequality (8) is invariant with respect to the change of
$\,\phi$ onto $1/\phi$ then one can assume further that both
$|w_j|\le 2$.

By Schottky' theorem (Cor.$\,$5 sec.2), $|f_j(z)|\le M(R,2),\,j=1,2$,
for $z\in B_R$.
By lemma 2 $\,|Vf|\le |V|_\kappa\cdot|f|(C+|\,log\,|f||)$ on $B_R$
for any vector field $V$ on $B$ with the constant
$C=2\,C_{\,0,1}+2\,log\,M,\,M:=M(R,3),$ if $w_1,w_2\in E\sm\{0,1\}$.
The same is true if one of $w_j$ is equal to 0 or 1 (see Prop.$\,$2).
Repeating literally the proof of Prop.$\,$2 (with this new constant $C$)
we obtain Harnack type inequalities
$$
\al\le\frac{C+|\,log\,|f(z)||}{C+|\,log\,|f(z_0)||}\le 1/\al\,,\
\quad z\in B_R\,.
$$

If $|w_1-w_2|>1$ then there is nothing to prove, the inequality (8) is valid
with any $C_R\ge 5$ and therefore we assume further that $|w_1-w_2|\le 1$
(and $|w_j|\le 2$).
If $|w_1-w_2|\ge e^{\,C}/(2\,e^{\,C}M)^{1/\al}$ then inequality (8)
is fulfilled with $C_R\ge 5^\al e^{\,C(1-\al)}2\,M$.
If $|w_1-w_2|<e^{\,C}/(2\,e^{\,C}M)^{1/\al}$ then it follows from the
inequality obtained above that $|f(z)|\le 1$
(otherwise $\al\cdot(C+log\,1/|w_1-w_2|)\le\,log\,e^{\,C}|f(z)|
\le\,log\,(2\,e^{\,C}M)$
in contradiction with the condition on $|w_1-w_2|$).
Therefore $\,\al\,log\,(e^{\,C}/|w_1-w_2|)\le\,log\,(e^{\,C}/|f(z)|)$, hence
$$
\frac{|(f_1-f_2)(z)|}{\sqrt{1+|f_1(z)|^2}\sqrt{1+|f_2(z)|^2}}<|f(z)|\le
e^{\,C(1-\al)}|w_1-w_2|\,^\al
$$
$$
\le 5^\al e^{\,C(1-\al)}
\left(\frac{|w_1-w_2|}{\sqrt{1+|w_1|^2}\sqrt{1+|w_2|^2}}\right)^\al\,
$$
and the inequality (8) is proved.

At last, if $|w_j|\le R\,'$ then, due to the Cor.$\,$4 sec.$\,$2,
$|f_j(z)|\le M(R,R\,')$
for $z\in B_R$, hence
$$
|(f_1-f_2)(z)|\le C_R\,(1+M(R,R\,')^2)\,|w_1-w_2|^{\,\al}\,.
$$
As the choice of $z_0$ on $B$ is arbitrary and $z_0\in B_R(z)$ if
$z\in B_R$ then, as it is proved,
$$
|w_1-w_2|\le C(R,R\,')\,|(f_1-f_2)(z)|^{\,\al}\,,
$$
and this is the left inequality in (9).
$\trr$
\vskip3mm

Lemma 2 together with Prop.$\,$4 show that for any normalized
holomorphic motion $\phi:B\times E\to\CC$
and any vector field $V$ bounded on $B_R$ the function
$V\phi$ satisfies on $B_R$ H\"older condition in $w$ with any exponent
$\beta<\al=e^{-R}$.
\vskip5mm

{\bf 4.$\,$Other generalizations.}
Harnack type inequalities are valid not only for holomorphic functions but
also for solutions of other elliptic equations and inequalities.
Let us start with evident generalizations for solutions of Beltrami equation.

Let a function $f$ in $\DD$ be a (generalized) solution of the equation
$f_{\bar z}=\mu f_z$ where $\mu\in L^\infty(\DD),\,\|\mu\|_\infty<1$.
Then there exists a quasiconformal homeomorphism
$\psi:\DD_z\to\DD_\zeta,\,\psi(0)=0,$ such that
$\psi_{\bar z}=\mu\,\psi_z$ in sense of distributions (see e.g. [A2]).
The function $h(\zeta)=f(\psi^{-1}(\zeta))$ is holomorphic in
$\DD,\,h(0)=f(0)$.
If $0<|f|<M$ then by the corollary in sec.1
$|f(z)|\le |f(0)|^{\,\beta(z)}\,M^{1-\beta(z)}$ where
$\beta(z)=\frac{1-|\psi(z)|}{1+|\psi(z)|}$.
Whereas if $f$ does not assume in $\DD$ the values 0 and 1
then by Prop.$\,$3
$$
\frac{1-|\psi(z)|}{1+|\psi(z)|}\,\le\,\frac{C_{\,0,1}+|\,log\,f(z)|}
{C_{\,0,1}+|\,log\,f(0)|}\,\le\,\frac{1+|\psi(z)|}{1-|\psi(z)|}\,,
$$
and corresponding inequalities for $|\,log|f(z)||$.

But these are evident generalizations.
More substantial is the following statement related with holomorphic
motions with a nonstandard complex structure in $B\times\CC$.
\vskip2mm

{\bf Proposition 5.} {\it Let $f$ be a continuous function without zeros
in the disk $\DD$ such that $|f|<1/e$ and the partial derivative
of $f_{\bar z}$ by $\bar z$ in sense of distributions is locally integrable
and satisfies almost everywhere in $\DD$ the inequality
\vskip1mm

\centerline{$|f_{\bar z}|\le A\,|f|\,log\,1/|f|$}
\vskip1mm

\np where $A$ is a function from $L^p(\DD)$ with some $p>2$.
Then there exists a constant $\,c\,$ depending only on $\,p$ and $\|A\|_{L^p}$
such that, for $z\in\DD\,$,}
\begin{equation}
e^{-c/(1-|z|)^{2-2/p}}\,\le
\frac{log\,1/|f(z)|}{log\,1/|f(0)|}\,\le\,e^{\,c/(1-|z|)}\,.
\label{10}
\end{equation}
\vskip1mm

$\trl$
One can assume that $f(0)>0$.

Denote by $log\,1/f$ the continuous logarithm of the function $1/f$ defined
by the condition $log\,1/f(0)>0$ and by $g$ denote the continuous logarithm
of the function $log\,1/f$ also defined by the condition $g(0)>0$.
As $Re\,log\,1/f>1$ then $|Im\,g|<\pi/2$.
As $log\,1/|f|\le|log\,1/f|$ then the function $g$ satisfies the inequality
$|g_{\bar\zeta}|\le A$, hence $g_{\bar\zeta}\in L^p(\DD)$.

The function
$$
a(z)=-\frac{1}{\pi}
\int_\DD\left(\,\frac{1}{\zeta-z}-\frac{1}{\zeta}\right)
g_{\bar\zeta}(\zeta)\,dS_\zeta
$$
is continuous on whole the plane (because $g_{\bar\zeta}\in L^p$), is bounded,
$|a|\le c_1(p)\|A\|_{L^p}$, $\,a(0)=0$ and the generalized derivative of
$a$ by $\bar z$ is equal to $g_{\bar z}$ (see [V] Ch.1).

By Weil lemma the function $g-a-g(0)=:h$ is holomorphic in $\DD$ and equal to
$0$ at $0$.
As $|Im\,g|<\pi/2$ and $|a|<c_0(p)\|A\|_p$ then
$|Im\,h|<c_0(p)\|A\|_p+\pi/2$ and therefore the function
$h$ has angular limit values almost everywhere on $\pa\DD$.
By Schwarz formuls
$$
  h(z)=\frac{i}{2\pi}\int_{\pa\DD}Im\,h(\zeta)\cdot\frac{\zeta+z}{\zeta-z}\,
|d\zeta|
$$
what follows that $\,|Re\,h|\le (\pi +2c_0(p)\|A\|_p)/(1-|z|)$, and the right
inequality (10) follows from the equality
 $|\,log\,1/f|=e^{\,Re(a+h)}\,log\,1/|f(0)|$.

The left inequality (10) is obtained from the right one using
automorphisms of $\DD$.
Let us fix an arbitrary point $z'\in\DD$ and set
$z=\psi(\zeta):=\frac{\zeta+z'}{1+\bar z'\zeta}$.
As $(f\circ\psi)_{\bar\zeta}=
f_{\bar z}\circ\psi\cdot\frac{1-|z'|^2}{(1-\bar z'\zeta)^2}\,$ then
$$
|(f\circ\psi)_{\bar\zeta}|\le A'|f\circ\psi|\,log\,1/|f\circ\psi|\quad
\text{with}\ \,
A'(\zeta)=\frac{1-|z'|^2}{|1+\bar z'\zeta|^2}A(\psi(\zeta))\in L^p(\DD)\,.
$$
As $\zeta=\frac{z-z'}{1-\bar z'z}$ then
$$
\|A'\|^p_{L^p}=
%\int_\DD (\frac{1-|z'|^2}{|1-z'\bar\zeta|^2}\,
%|A(\psi(\zeta)|)^p\,dS_\zeta=
\int_\DD \frac{|1-\bar z'z|^{2p-4}}{(1-|z'|^2)^{p-2}}\,|A(z)|^p\,dS_z
\,\le\,\frac{2^{\,p-2}\|A\|^p_{L^p}}{(1-|z'|)^{p-2}}\,.
$$

Let $\al=arg\,f(z')\in(-\pi,\pi]$ and
$log\,e^{i\al}/f,\ log\,log\,e^{i\al}/f$ are continuous branches of
logarithms which are positive at $\zeta=0$ (corresponding to $z=z'$).
As it is proved above
$$
\frac{log\,1/|f\circ\psi(\zeta)|}{log\,1/|f\circ\psi(0)|}\,\le\,
exp\,[(\pi+3c_0(p)\|A'\|_{L^p})/(1-|\zeta|)]
$$
\np As $\psi(0)=z'$ and $\psi(-z')=0$ it follows at $\zeta=-z'$ that
$$
log\,\frac{1}{|f(0)|}\,\le\,e^{\,c/(1-|z'|)^{\,2-2/p}}\,log\,\frac{1}{|f(z')|}
$$
\np and the left inequality (10) is proved.
$\trr$
\vskip2mm

More symmetric inequalities would be obtained by $p=2$ but $c_0(p)\to\infty$
as $\,p\downarrow 2$.
\vskip3mm

Similar estimates are valid also for functions on the whole plane but they
depend not on geometry (Kobayashi distance on the plane vanishes identically)
but on estimates of $f_{\bar z}$.
\vskip3mm

{\bf Proposition 6.} {\it Let $f$ be a continuous function without zeros on
the plane $\CC_z$ such that $|f|<1$ and the partial derivative $f_{\bar z}$
by $\bar z$ in sense of distributions is locally integrable and satisfies
almost everywhere on $\CC$ the inequality
\vskip1mm

\centerline{$|f_{\bar z}|\le A\,|f|\cdot|\,log\,M/f|$}
\vskip1mm

\np where $A$ is a function from $(L^p\cap L^{p'})(\CC)$
with some $p>2,\,p'<2$, $M\ge 1$ is a constant and
$\,log\,M/f$ is a continuous branch of logarithm.
Then $f$ extends continuously onto $\hat\CC$ and there exists a constant
$C>0$ depending only on $\,p, p'$ and
$\|A\|_{L^p}+\|A\|_{L^{p'}}$ such that}
\begin{equation}
C^{-1}\le\,\frac{|\,log\,M/f(z)|}{|\,log\,M/f(0)|}\,\le\,C\,.
\label{11}
\end{equation}

$\trl$
Set $g=\,log\,M/f$.
From the estimate of $f_{\bar z}$ we obtain that
$\,|g_{\bar z}|\le A\,|g|$ and $\,Re\,g>0$.

Set $a=g_{\bar z}/g$.
Then $a\in (L^p\cap L^{p'})(\CC)$ and $g_{\bar z}=a\,g$.
The function
$$
\hat a(z)=-\frac{1}{\pi}
\int_\CC\left(\frac{1}{\zeta-z}-\frac{1}{\zeta}\right)
a(\zeta)\,dS_\zeta
$$
is continuous (as $a\in L^p$) and bounded
(as $a\in L^{p'}$) on whole the plane $\CC$ and the generalized derivative
of $\hat a$ by $\bar z$ equals $a$ (see [V] Ch.1), besides
$a(0)=0$ and $\,|\hat a|\le C_1\,(\|A\|_{L^p}+\|A\|_{L^{p'}})$
with a constant $C_1$.

For generalized derivative by $\bar z$ of the function $ge^{-\hat a}$
the Leibnitz rule is applicable:
if $g^\eps,a^\eps$ is a standard mollification with a smooth ``hat'' then
$g^\eps e^{-\hat{a^\eps}}\to ge^{-\hat a}$ in sense of distributions, hence
$$
(ge^{-\hat a})_{\bar z}=
lim_{\eps\to 0}(g^\eps_{\bar z}e^{-\hat{a^\eps}}-
g^\eps a^\eps e^{-\hat{a^\eps}})=(g_{\bar z}-a\,g)e^{-\hat a}\,.
$$

\np As $\,g_{\bar z}=a\,g\,$ it follows from this and Weil lemma that
$ge^{-\hat a}$ is a holomorphic function on whole the plane and it does not
have zeros.
And then the function $log\,g-\hat a=:h\,$ is also holomorphic on whole the
plane, has uniformly bounded imaginary part, hence is constant.
As $\hat a(0)=0$ then $h\equiv log\,g(0)$, hence
$g=g(0)\,e^{\hat a}$.
\vskip1mm

It follows that $\,log\,M/ff=e^{\hat a}\,log\,M/f(0)\,$
and the right inequality
(11) is fulfilled with the constant $C=sup\,e^{\,Re\,\hat a}$.
As the conditions of the proposition do not change by shifts
$z\mapsto z+const$ then the left inequality follows also.
$\trr$
\vskip2mm

{\bf Corollary 1.} {\it If $\ \underline{lim}_{\,z\to\infty}\,|f|=0$
then $\,f\equiv 0$.}
\vskip1mm

$\trl$
Apply the right inequality (11) to $f(z+R)$ with arbitrary big $R$.
$\trr$
\vskip2mm

{\bf Corollary 2.} $\quad sup\,|f|\,\le\,(inf\,|f|)^{1/C}\,.$
\vskip1mm

$\trl$
We can assume that $\,inf\,|f|=f(0)>0$ (shift by $z$ and multiply $f$ on
 a constant equal to one by modulus).
$\trr$
\vskip2mm

{\bf Corollary 3.} {\it If $\,|f_{\bar z}|\le A\,|f|\,$ with a function
$A$ as in proposition $6$ then
$$
sup\,|f|\,\le\,C(A)\cdot inf\,|f|
$$
with similar constant $C(A)>1$.}
\vskip1mm

$\trl$
It follows from given Lipschitz condition that
 $f/e$ satisfies the log-Lipschitz condition from the proposition and
therefore
$sup\,|f|^{1/C}\le e^{1-1/C}\cdot inf\,|f|^{1/C}$.
$\trr$
\vskip3mm

Propositions 5, 6 and corollaries are valid surely
also for functions satisfying the inequalities from these
propositions in which the Cauchy -- Riemann operator $\,\pa/\pa\bar z$
is substituted by the Beltrami operator $\pa/\pa\bar z-\mu\,\pa/\pa z$
(the proofs are reduced to Prop.$\,$5, 6 by evident quasicon\-formal
transformations of $\DD$ or $\CC$).
\vskip10mm

\centerline{* * * * *}
\vskip4mm

\centerline{\large Referencies}
\bigskip

\np [Ag]
$\ $
S.Agard, Distorsion theorems for quasiconformal mappings,
Ann. Acad. Sci. Fenn. AI, ${\bf 413}$, 1968, 1--7.
\smallskip

\np
[A1]
$\ $
L.Ahlfors, Conformal Invariants, Topics in Geometric Function Theory,
AMS Chelsea Publ., Providence, RI, 2010.
\smallskip

\np
[A2]
$\ $
L.Ahlfors, Lectures on Quasiconformal Mappings, Van Nostrand,
Prin\-ceton, NJ, 1966.
\smallskip

\np
[GJW]
F.Gardiner, Y.Jiang, Z.Wang, Holomorphic motions and related topics,
Geometry of Riemann surfaces, London Math. Soc. Lecture Notes Ser.,
368, Cambridge Univ. Press, Cambridge, 2010, 153--193
(arXiv:mathCV 0802.2111, 2008, 32pp.).
\smallskip

\np
[HK]
W.K.Hayman, P.B.Kennedy, Subharmonic functions, vol.1, Academic Press,
London, 1976.

\np
[H]$\,$ J.A.Hempel, The Poincar\'e metric on the twice punctured plane
and the theorems of Lamdau and Schottky, J. London Math. Soc. (2),
{\bf 20}, 1979, 435--445.
\smallskip

\np
[J]$\,$ J.A.Jenkins, On explicit bounds in Schottky's theorem,
Canad. J. Math., {\bf 7}, 1955, 76--82.
\smallskip

\np
[Ko]
$\ $
S.Kobayashi, Hyperbolic Manifolds and Holomorphic Mappings,
\newline M.Decker, NY, 1970.
\smallskip

\np [LV]$\,$
O.Lehto, K.I.Virtanen, On the existence of quasiconformal mappings with
prescribed complex dilatation, Ann. Acad. Sci. Fenn. AI, {\bf 274},
1960, 1--24.
\smallskip

\np
[MSS]
$\ $
R.Ma\~n\'e, P.Sad, D.Sullivan, On the dynamics of rational maps,
Ann. Sci. ENS, {\bf 16}, 1983, 193--217.
\smallskip

\np [N]$\ $ Z.Nehari, Conformal Mappings, Dover Publ., NY, 1975.
\smallskip

\np
[R1]
$\ $
H.Royden, Remarks on Kobayashi metric, Proc. Maryland Conference on Several
Complex Variables, Springer Lecture Notes, {\bf 185}, Springer-Verlag, Berlin,
1971.
\smallskip

\np
[R2]
$\ $
-- -- , The extension of regular holomorphic maps,
Proc.Amer.\-Math.Soc., {\bf 43}, 1974, 306--310.
\smallskip

\np
[V]
$\ $
I.N.Vekua, Generalized analytic functions, Pergamon Press, Oxford, 1962.
Russian original: ``Nauka'', Moscow, 1988.
\vskip10mm

\np chirka@mi.ras.ru
\vskip1mm

\np Steklov Math. Institute

\np Moscow

\end{document}